\documentclass[10pt]{article}
\usepackage{amsmath,amssymb}
\newcommand{\dis}{\displaystyle}
\newcommand{\pf}{\noindent{\bf Proof. }}
\newcommand{\nint}{\int_{\R^N}}
\newcommand{\eps}{{\varepsilon}}
\newcommand{\Div}{\mathrm{div}}
\newcommand{\dt}{{\partial_t}}
\newcommand{\N}{\mathbb N}
\newcommand{\R}{\mathbb{R}}

\newcommand{\te}{{\theta}}
\newcommand{\tek}{{\theta_k}}
\newcommand{\tes}{\theta^*}
\newcommand{\tesk}{{\theta_k^*}}

\def\na{\nabla}
\def\pa{\partial}
\newcommand{\tP}{\tilde{P}}
\newcommand{\La}{{\Lambda}}
\newcommand{\un}{{\mathbf{1}}}
\newcommand{\Cb}{{\overline{C}}}
\newcommand{\Qks}{{\overline{B}_k^*}}
\newcommand{\At}{{\mathcal{A}}(t)}
\newcommand{\Bt}{{\mathcal{B}}(t)}
\newcommand{\Ct}{{\mathcal{C}}(t)}
\newcommand{\tbk}{\overline{\te}_k}
\newcommand{\tbku}{\overline{\te}_{k-1}}
\newcommand{\tbsk}{\overline{\te}^*_k}
\newcommand{\tbsku}{\overline{\te}^*_{k-1}}
\newcommand{\ds}{\displaystyle}

\def\qed{\hbox{${\vcenter{\vbox{
  \hrule height 0.4pt\hbox{\vrule width 0.4pt height 6pt
  \kern5pt\vrule width 0.4pt}\hrule height 0.4pt}}}$}}

\newtheorem{theo}{Theorem}
\newtheorem{prop}[theo]{Proposition}
\newtheorem{lemm}[theo]{Lemma}
\newtheorem{coro}[theo]{Corollary}

\author{Luis Caffarelli\thanks{caffarel@math.utexas.edu},
Alexis Vasseur\thanks{vasseur@math.utexas.edu}\\
 \\
\small Department of Mathematics\\
\small University of Texas at Austin
 }
\date{}
\title{Drift diffusion equations with fractional diffusion
and the quasi-geostrophic equation}

\begin{document}
\maketitle

\bibliographystyle{plain}

 {\small \noindent{\bf Abstract:} Motivated by the critical
 dissipative quasi-geostrophic equation, we prove that
 drift-diffusion equations with $L^2$ initial data and minimal
 assumptions on the drift are locally Holder continuous. As an
 application we show that solutions of the quasi-geostrophic
 equation with initial $L^2$ data and critical diffusion
 $(-\Delta)^{1/2}$, are locally smooth for any space dimension.
 }

\section{Introduction}

Non linear evolution equations with fractional diffusion arise in
many contexts: In the quasi-geostrophic flow model (Constantin
\cite{Constantin_Lecture}), in boundary control problems
(Duvaut-Lions \cite{Duvaut_Lions}), in surface flame propagation and
in financial mathematics. In this paper, motivated by the
quasi-geostrophic model, we study the equation:
\begin{equation}\label{eq_quasi_linear}
\begin{array}{l}
\ds{\dt\te+v\cdot\na\te=-\La\te,\qquad x\in\R^N,}\\[0.3cm]
\ds{\Div v=0,}
\end{array}
\end{equation}
where $\La\te=(-\Delta)^{1/2}\te$.
The main two theorems are roughly the following a priori estimates:
\begin{theo}\label{theo L^2 L^infty}
(From $L^2$ to $L^\infty$). Let $\te(t,x)$ be a function in
$L^\infty(0,T;L^2(\R^N))\cap L^2(0,T;H^{1/2}(\R^N))$. For every
$\lambda>0$ we define:
$$
\te_\lambda=(\te-\lambda)_+.
$$
If $\te$ (and $-\te$) satisfies for every $\lambda>0$ the level set
energy inequalities:
\begin{eqnarray*}
&&\nint\te_\lambda^2(t_2,x)\,dx+2\int_{t_1}^{t_2}\nint|\La^{1/2}\te_\lambda|^2\,dx\,dt\\
&&\qquad\leq\nint\te_\lambda^2(t_1,x)\,dx,\qquad 0<t_1<t_2,
\end{eqnarray*}
then:
$$
\sup_{x\in\R^N}|\te(T,x)|\leq C^*\frac{\|\te_0\|_{L^2}}{T^{N/2}}.
$$
\end{theo}
\noindent{\bf Remark:} That solutions to equation
(\ref{eq_quasi_linear}) are expected to satisfy the energy
inequality follows from writing $\La$ as the normal derivative of
the harmonic extension of $\te$ to the upper half space. Existence
theory is sketching in appendix C. 
In the case of Quasi-geostrophic equation it can also be seen as a
corollary of Cordoba and Cordoba \cite{Cordoba^2_2004}.

Those energy inequalities are reminiscent of the notion of entropic
solutions for scalar conservation laws.  Consider a weak solution of
(\ref{eq_quasi_linear}) lying in $L^2(H^{1/2})$ and for which we can
define the equality (in the sense of distribution for example):
$$
\phi'(\te)v\cdot\nabla\te=\Div(v\phi(\te)),
$$
for any Lipschitz function $\phi$. Then $\te$ verifies the level set
energy inequalities. In the case of the Quasi-geostrophic equation,
$v\in L^2(H^{1/2})$ and we can give a sense to:
$$
v\cdot\nabla\phi(\te).
$$
Indeed, using the harmonic extension, it can be shown that if $\te$
lies in  $L^2(H^{1/2})$ so does $\phi(\te)$. and so
$\nabla\phi(\te)$ lies in $L^2(H^{-1/2})$.
 \vskip0.3cm

For the second theorem, (from $L^\infty$ to $C^\alpha$), we need
better control of $v$:

\begin{theo}\label{theo L^infty to C^alpha}
(From $L^\infty$ to $C^\alpha$). We define $Q_r=[-r,0]\times
[-r,r]^N$, for $r>0$. Assume now that $\te(t,x)$ is bounded in
$[-1,0]\times\R^N$ and $v|_{Q_1}\in L^\infty(-1,0;BMO)$, then $\te$
is $C^\alpha$ in $Q_{1/2}$.
\end{theo}
\noindent{\bf Remark1:} The global bound of $\te$ is not really
necessary, only local $L^\infty$ and integrability at infinity
against the Poisson kernel, as we will see later. \vskip0.2cm

\noindent{\bf Remark2:} Note that both theorems depend only on the
resulting energy inequality and not on the special form of $\La$.
\vskip0.2cm

From these two theorems, the regularity of solutions to the
quasi-geostrophic equation follows.
\begin{theo}\label{theo quasi}
Let $\te$ be a solution to an equation
\begin{equation}\label{eq:quasi}
\begin{array}{l}
\ds{\dt\te+u\cdot\na\te=-\La\te,\qquad x\in\R^N,}\\[0.3cm]
\ds{\Div u=0,}
 \end{array}
\end{equation}
with
 \begin{equation}\label{eq:uquasi}
 u_j=\overline{R}_j[\te],
\end{equation}
  $\overline{R}_j$ a singular integral
operator. Assume also that $\te$ verifies the level set energy
inequalities stated in Theorem \ref{theo L^2 L^infty}. Then, for
every $t_0>0$ there exists $\alpha$ such that $\te$ is bounded in
$C^\alpha([t_0,\infty[\times\R^N)$.
\end{theo}
Indeed, Theorem \ref{theo L^2 L^infty} gives that $\te$ is uniformly
bounded  on $[t_0,\infty[$ for every $t_0>0$. Singular integral
operators are bounded from $L^\infty$ to $BMO$. This gives that
$u\in L^\infty(t_0,\infty; BMO(\R^N))$ and, after proper scaling,
Theorem \ref{theo L^infty to C^alpha} gives the result of Theorem
\ref{theo quasi}. \vskip0.2cm

\noindent{\bf Remark 1:} Higher regularity then follows from
standard potential theory, by noticing that the fundamental solution
of the operator:
$$
\dt+\La\te=0
$$
is the Poisson kernel and that in the non linear term we can
substract a constant both $\te_0$ from $\te$ and $u_0$ from $u$,
this last one by a change of coordinates:
$$
x^*=x-tu_0,
$$
doubling its Holder decay (see appendix). \vskip0.2cm

Strictly speaking, the dissipative quasi-geostrophic flow model in
the critical case corresponds to the case $N=2$ and
\begin{equation*}
\begin{array}{l}
\dis{u_1=-R_2\te,}\\
\dis{u_2=R_1\te},
\end{array}
\end{equation*}
where $R_i$ is the usual Riesz transform defined from the Fourier
transform: $\widehat{R_i\te}=\frac{i\xi_i}{|\xi|}\widehat{\te}$.
This model was introduced by some authors as a toy model to
investigate the global regularity of solutions to 3D fluid mechanics
(see for instance \cite{Constantin_Lecture}). When replacing the
diffusion term $-\La$ by $-\La^\beta$, $0\leq\beta\leq2$, the
situation is classically decomposed into 3 cases: The subcritical
case for $\beta>1$, the critical case for $\beta=1$ and the
supercritical case for $\beta<1$.

Weak solutions has been constructed by  Resnick  in \cite{Resnick}.
Constantin and Wu showed in \cite{Constantin_Wu} that in the
subcritical case any solution with smooth initial value is smooth
for all time. Constantin Cordoba and Wu showed in
\cite{Constantin_Cordoba_Wu} that the regularity is conserved for
all time in the critical case provided that the initial value is
small in $L^\infty$. In both  the critical case and supercritical
cases, Chae and Lee considered in \cite{Chae_Lee} the well-posedness
of solutions with initial conditions small in Besov spaces (see also
Wu \cite{Wu_Besov_small}).

Notice that our case corresponds to the critical case and global
regularity in $C^{1,\beta}$, $\beta<1$ is showed for any initial
value in the energy space without hypothesis of smallness. This
ensures that the solutions are classical.

Let us also cite a result of maximum principle due to Cordoba and
Cordoba \cite{Cordoba^2_2004},  results of behavior in large time
due to Schonbeck and Schonbeck \cite{Schonbek^1_2003},
\cite{Schonbek^2}, and a criteria for blow-up in Chae
\cite{Chae2006}.

\noindent{\bf Remark 2:} In a recently posted preprint in arXiv,
Kiselev, Nazarov, and Volberg present a very elegant proof of the
fact that in 2D, solutions with periodic $C^\infty$ data for the
quasi-geostrophic equation remain $C^\infty$ for all time
(\cite{Kislev}).\vskip0.2cm

We conclude our introduction by pointing out that our techniques
also can be seen as a parabolic De Giorgi Nash Moser method to treat
"boundary parabolic problems" of the type:
\begin{eqnarray*}
&&\Div (a\nabla \te)=0,\qquad \mathrm{in}\ \ \Omega\times[0,T]\\
&&[f(\te)]_t=\te_\nu\qquad \mathrm{on} \ \
\partial\Omega\times[0,T],
\end{eqnarray*}
 that
arise in boundary control (see Duvaut Lions \cite{Duvaut_Lions}).
Note also that similar results to Theorem \ref{theo L^2 L^infty} can
be obtained even for systems (See Vasseur \cite{Vasseur} and Mellet,
Vasseur \cite{Mellet} for  applications of the method in fluid
mechanics).

\section{$L^\infty$ bounds}

This section is devoted to the proof of Theorem \ref{theo L^2
L^infty}. We use the level set energy inequality  for:
$$
\lambda=C_k=M(1-2^{-k}),
$$
where $M$ will be chosen later. This leads to the following energy
inequality for the level set function $\tek=(\te-C_k)_+$:
\begin{equation}\label{eq_tek}
\dt\nint\tek^2\,dx+2\nint|\La^{1/2}\tek|^2\,dx\leq0.
\end{equation}
Let us fix a $t_0>0$, we want to show that $\te$ is bounded for
$t>t_0$. We introduce $T_k=t_0(1-2^{-k})$, and the level set of
energy/dissipation of energy:
$$
U_k=\sup_{t\geq
T_k}\left(\nint\tek^2\,dx\right)+2\int_{T_k}^\infty\nint|\La^{1/2}\tek|^2\,dx\,dt.
$$
integrating (\ref{eq_tek}) in time between $s$, $T_{k-1}<s<T_k$, and
$t>T_k$ and between $s$ and $+\infty$ we find:
$$
U_k\leq 2\nint\tek^2(s)\,dx.
$$
Taking the mean value in $s$ on $[T_{k-1},T_k]$ we find:
\begin{equation}\label{eq_Uk}
U_k\leq
\frac{2^{k+1}}{t_0}\int_{T_{k-1}}^{\infty}\nint\tek^2\,dx\,dt.
\end{equation}
We want to control the right-hand side by $U_{k-1}$ in a non linear
way.  Sobolev and Holder inequalities give:
$$
U_{k-1}\geq C
\|\te_{k-1}\|^2_{L^{\frac{2(N+1)}{N}}(]T_{k-1},\infty[\times\R^N)}.
$$
Note that  if $\tek>0$ then $\te_{k-1}\geq 2^{-k}M$. So
$$
\un_{\{\tek>0\}}\leq \left(\frac{2^k}{M}\te_{k-1}\right)^{2/N}.
$$
Hence:
\begin{eqnarray*}
&&U_k\leq
\frac{2^{k+1}}{t_0}\int_{T_{k-1}}^{\infty}\nint\te_{k-1}^2\un_{\{\tek>0\}}\,dx\,dt\\
&&\qquad\leq
2\frac{2^{\frac{N+2}{N}k}}{t_0M^{2/N}}\int_{T_{k-1}}^{\infty}\nint\te_{k-1}^{2\frac{N+1}{N}}\,dx\,dt\\
&&\qquad\leq
2C\frac{2^{\frac{N+2}{N}k}}{t_0M^{2/N}}U_{k-1}^{\frac{N+1}{N}}.
\end{eqnarray*}
For $M$ such that $M/t_0^{N/2}$ is big enough (depending on  $U_0$)
we have $U_k$ which converges to 0. This gives $\te\leq M$ for
$t\geq t_0$. The same proof on $-\te$ gives the same bound
for$|\te|$. Note that $U_0\leq\|\te_0\|^2_{L^2}$. The scaling
invariance $\te_\eps(s,y)=\te(\eps s,\eps y)$ gives the final
dependence with respect to $\|\te_0\|_{L^2}$.\qquad\qed \vskip0.2cm

This theorem leads to the following corollary.
\begin{coro}\label{lemm_linfty}
There exists a constant $C^*>0$ such that  any solution $\te$ of
(\ref{eq:quasi}) (\ref{eq:uquasi})  verifies:
\begin{eqnarray*}
&& \sup_{x\in\R^N}|\te(T,x)|\leq
C^*\frac{\|\te_0\|_{L^2(\R^N)}}{T^{N/2}},\\
&&\|u(T,\cdot)\|_{BMO(\R^N)}\leq
C^*\frac{\|\te_0\|_{L^2(\R^N)}}{T^{N/2}}.
\end{eqnarray*}
\end{coro}

\pf First note that the property on $u$ follows directly from the
property on $\te$ and the imbedding of the Riesz function from
$L^\infty$ to $BMO$. We make use of the following result of Cordoba
and Cordoba (see \cite{Cordoba^2_2004}): for any convex function
$\phi$ we have the pointwise inequality:
$$
-\phi'(\theta)\La\te\leq-\La(\phi(\te)).
$$
Making use of this inequality with:
$$
\phi_k(\te)=(\te-C_k)_+=\tek
$$
leads to:
$$
\dt \tek+u\cdot\nabla\tek\leq-\La\tek.
$$
Multiplying by $\tek$ and integrating in $x$ gives (\ref{eq_tek}),
using that $u$ is divergence free. \qquad\qed \vskip0.3cm

\noindent{\bf Remark:} We point out that the level set energy
inequalities we assume in Theorem \ref{theo L^2 L^infty} is
heuristically a general fact (See appendix C).

\section{Local energy inequality:}

To get the $C^\alpha$ regularity we need to get an energy inequality
which is local in time and space. Due to the non locality of the
diffusion operator this cannot be obtained directly. It relies on
the following classical representation of the operator $\La$. We
introduce first the harmonic extension $L$ defined from
$C^{\infty}_{0}(\R^N)$ to $C^{\infty}_0(\R^N\times\R^+)$ by:
\begin{eqnarray*}
&&-\Delta L(\te)=0 \qquad \mathrm{in} \ \ \R^N\times(0,\infty),\\
&&L(\te)(x,0)=\te(x) \qquad\mathrm{for} \ \ x\in\R^N.
\end{eqnarray*}
Then the following result holds true:
consider  $\te$ defined on $\R^N$. then:
\begin{equation}\label{lemm:La}
\La\te(x)=\pa_\nu[L\te](x),
\end{equation}
where we denote $\pa_\nu[L\te]$ the normal derivative of $L\te$ on
the boundary $\{(x,0)| x\in\R^N\}$.

In the following, we will denote:
\begin{equation}\label{eq_def_tes}
\tes(t,x,z)=L(\te(t,\cdot))(x,z).
\end{equation}

%

We denote $B_r=[-r,r]^N$,
$B^*_r=B_r\times(0,r)\in\R^N\times(0,\infty)$, and
$[y]_+=\sup(0,y)$.

The rest of this section is devoted to the proof of the following
proposition:
\begin{prop}\label{prop_local_energy teC}
Let $t_1,t_2$ be such that $t_1<t_2$ and  let $\te\in
L^\infty(t_1,t_2;L^2(\R^N))$ with $\La^{1/2} \te\in
L^2((t_1,t_2)\times\R^N)$, be solution to (\ref{eq_quasi_linear})
with a velocity $v$ satisfying:
\begin{equation}\label{condition sur v}
\|v\|_{L^\infty(t_1,t_2;BMO(\R^N))}+\sup_{t_1\leq t\leq
t_2}\left|\int_{B_2}v(t,x)\,dx\right|\leq C.
\end{equation}
Then there exists a constant $\Phi$ (depending only on $C$) such
that for every $t_1\leq t\leq t_2$ and cut-off function $\eta$
compactly supported in $B_2^*$:
\begin{equation}\label{eq_energy_locale}
\begin{array}{l}
 \dis{\qquad \int_{t_1}^{t_2}\int_{
B^*_2}|\nabla(\eta[\tes]_+)|^2\,dx\,dz\,dt+\int_{B_2}(\eta[\te]_+)^2(t_2,x)\,dx}\\[0.3cm]
\dis{\leq
\int_{B_2}(\eta[\te]_+)^2(t_1,x)\,dx+\Phi\int_{t_1}^{t_2}\int_{B_2}([\na\eta][\te]_+)^2\,dx\,dt}\\[0.3cm]
\ds{\qquad\qquad+\int_{t_1}^{t_2}\int_{B^*_2}([\na\eta][\tes]_+)^2\,dx\,dz\,dt.}
\end{array}
\end{equation}
\end{prop}
\pf 
We have for every $t_1<t<t_2$:
\begin{eqnarray*}
0&=&\int_0^\infty\int_{\R^N} \eta^2[\tes]_+\Delta\tes\,dx\,dz\\
&=&-\int_0^\infty \int_{\R^N}|\nabla(
\eta[\tes]_+)|^2\,dx\,dz+\int_0^\infty \int_{\R^N}|\nabla\eta|^2
[\tes]_+^2\,dx\,dz\\
&&\qquad +\int_{\R^N}\eta^2[\te]_+\La\te\,dx.
\end{eqnarray*}
Using equation (\ref{eq_quasi_linear}), we find that:
\begin{eqnarray*}
&&\qquad-\int_{\R^N}\eta^2[\te]_+\La\te\,dx\\
&&=\frac{\partial}{\partial
t}\left(\int_{\R^N}\eta^2\frac{[\te]_+^2}{2}\,dx\right)-\int_{\R^N}\nabla{\eta^2}\cdot
v\frac{[\te]_+^2}{2}\,dx.
\end{eqnarray*}
This leads to:
\begin{eqnarray*}
&&\int_{t_1}^{t_2}\int_0^\infty\int_{\R^N} |\nabla(
\eta[\tes]_+)|^2\,dx\,dz\,ds+\int_{\R^N}\eta^2\frac{[\te]_+^2(t_2)}{2}\,dx\\
&&\leq
\int_{\R^N}\eta^2\frac{[\te]_+^2(t_1)}{2}\,dx+\int_{t_1}^{t_2}\int_0^\infty\nint
|\nabla\eta|^2
[\tes]_+^2\,dx\,dz\,ds\\
&&\qquad+\left|\int_{t_1}^{t_2}\int_{\R^N}\eta\nabla{\eta}\cdot
v[\te]_+^2\,dx\,ds\right|.
\end{eqnarray*}
To dominate the last term, we first use the Trace theorem and
Sobolev imbedding to find:
\begin{eqnarray*}
\|\eta \te_+\|^2_{L^{\frac{2N}{N-1}}(\R^N)}&\leq C&\|\eta
\te_+\|^2_{H^{1/2}(\R^N)}=C\int_{B_2}(\eta\tes_+)\La(\eta\tes_+)\,dx\\
&=& C\int_{B_2^*}|\na L(\eta(\tes)_+)|^2\,dx\,dz\\
 &\leq&C\int_{B_2^*}|\na[\eta(\tes)_+]|^2\,dx\,dz.
\end{eqnarray*}
In the last inequality we have used the fact that $L(\eta(\tes)_+)$
is harmonic and have the same trace than $\eta(\tes)_+$ at $z=0$.
Therefore we split:
\begin{eqnarray*}
&&\qquad\left|\int_{t_1}^{t_2}\int_{\R^N}\nabla{\eta^2}\cdot
v\frac{[\te]_+^2}{2}\,dx\,ds\right|\\
&&\leq \varepsilon \int_{t_1}^{t_2}\|\eta
\te_+\|^2_{L^{\frac{2N}{N-1}}(\R^N)}\,ds+\frac{1}{\varepsilon}\int_{t_1}^{t_2}\|[\na\eta]
v[\te]_+\|^2_{L^{\frac{2N}{N+1}}}\,ds.
\end{eqnarray*}
The first term is absorbed by the left. The second can be bounded,
using Holder inequality, by:
$$
\frac{1}{\varepsilon}\|v\|^2_{L^{\infty}(t_1,t_2;L^{2N}(B_2))}
\int_{t_1}^{t_2}\int_{\R^N}|[\na\eta][\te]_+|^2\,dx\,ds,
$$
which gives the desired result. \qquad\qed

\section{From $L^2$ to $L^\infty$:}\label{section_lemme1}

In these two sections (\ref{section_lemme1} and
\ref{section_lemme2}) we follow De Giorgi's ideas in his
"oscillation lemma" (see \cite{DeGiorgi}) to prove Holder
continuity: Suppose that $\te$ oscillates in $Q_1=[-1,0]\times B_1$
between $-2$ and 2, but it is negative most of the time. In
particular, if $\|\te_+\|_{L^2}$ is very small, then we prove that
$\te_+|_{Q_{1/2}=[-1/2,0]\times B_{1/2}}\leq 2-\lambda$, effectively
reducing the oscillations of $\theta$ by $\lambda$ (see section
\ref{section_lemme1}). Of course, we do not know a priori that this
is the case. But we do know that in $Q_1$, $\te$ is at least half of
the time positive, or negative, say negative. Then we reproduce a
version of De Giorgi's isoperimetric inequality that says that to go
from zero to one $\te$ needs "some room" (section
\ref{section_lemme2}). Therefore the set $\{\te\leq 1\}$ is
"strictly larger" than the set $\{\te\leq 0\}$. Repeating this
argument at truncation levels $\sigma_k=2-2^{-k}$, we fall, after a
finite number of steps, $k_0$, into the first case, effectively
diminishing the oscillations of $\te$ by $\lambda 2^{-k_0}$. This
implies Holder continuity (section \ref{section_oscillation lemma}
and section \ref{section_theorem}).

This section is devoted to the proof of the following technical
lemma. It says that, under suitable conditions on $v$, we can
control the $L^\infty$ norm of $\te$ from the $L^2$ norm of both
$\te$ and $\te^*$ locally.
\begin{lemm}\label{lemm_L^2_L^infty} There exists
$\eps_0>0$ and $\lambda>0$ such that for every $\te$ solution to
(\ref{eq_quasi_linear}) with a velocity $v$ satisfying:
\begin{eqnarray*}
\|v\|_{L^\infty(-4,0;BMO(\R^N))}+\sup_{-4\leq t\leq
0}\left|\int_{B_4}v(t,x)\,dx\right|\leq C,
\end{eqnarray*}
the following property holds true.

If we have:
$$
\tes\leq2\qquad \mathrm{in} \ \ [-4,0]\times B^*_4,
$$
and
$$
\int_{-4}^0\int_{B_4^*}(\tes)^2_+\,dx\,dz\,ds+
\int_{-4}^0\int_{B_4}(\te)^2_+\,dx\,ds\leq \eps_0,
$$
then:
$$
(\te)_+\leq 2-\lambda \qquad \mathrm{on} \ \ [-1,0]\times B_1.
$$
\end{lemm}

\pf We split the proof of the lemma into several steps. \vskip0.2cm
\noindent {\bf Step1. Useful barrier functions and setting of the
constant $\lambda$:} Consider the function $b_1$, defined by:
\begin{eqnarray*}
&&\Delta b_1=0\qquad\mathrm{in} \ \ B^*_4\\
&&b_1=2\qquad \mathrm{on \ the \ sides \ of \ the \ cube \ } B^*_4 \
\mathrm{ except \ for } \
z=0\\
&& b_1=0 \qquad \mathrm{for} \  z=0.
\end{eqnarray*}
Then there exists $\lambda>0$ such that:
$$
b_1(x,z)\leq 2-4\lambda \ \mathrm{on} \ B^*_2.
$$
This result follows directly from the maximum principle. We consider
now $b_2$ harmonic function defined by:
\begin{eqnarray*}
&&\Delta b_2=0\qquad\mathrm{in} \ \ [0,\infty[\times[0,1],\\
&&b_2(0,z)=2\qquad 0\leq z\leq 1,\\
&& b_2(x,0)=b_2(x,1)=0\qquad 0< x<\infty.
\end{eqnarray*}
Then there exists $\overline{C}>0$ such that:
\begin{equation}\label{eq_b_2}
|b_2(x,z)|\leq \overline{C}e^{-\pi x}.
\end{equation}
Notice that $\overline{C}$ is universal. Actually $b_2$ can be
explicitly computed by the method of separation of variables:
$$
b_2(x,z)=\sum_{p=0}^\infty\frac{4}{\pi(2p+1)}e^{-[(2p+1)\pi]x}\sin((2p+1)\pi
z).
$$
\vskip0.2cm\noindent{\bf Step 2. Setting of constants:} In this step
we fix a set of constants. We make the choice to set them right away
to convince the reader that there is no loop in the proof.

\begin{lemm}
There exist $0<\delta<1$ and $M>1$ such that for every $k>0$:
\begin{eqnarray*}
&&N\Cb e^{-\pi\frac{2^{-k}}{\delta^k}}\leq\lambda 2^{-k-2},\\
&&\frac{M^{-k}}{\delta^{k+1}}\|P(1)\|_{L^2}\leq\lambda 2^{-k-2},\\
&&M^{-k}\geq C_0^k M^{-(1+1/N)(k-3)}\qquad k\geq 12N,
\end{eqnarray*}
where $\Cb$ is defined from step 1, $P(1)$ is the the value at $z=1$
of the Poisson kernel $P(z)(x)$, and $C_0$ is defined by
(\ref{eq_C_0}).
\end{lemm}
The proof is easy. We construct first $\delta$ to verify the first
inequality in the following way. If $\delta< 1/4$, the inequality is
true for $k>k_0$ due to the exponential decay. If necessary, we then choose $\delta$ smaller to make the inequality also valid for $k<k_0$. Now that  $\delta$
has been fixed, we have to choose  $M$ large  to satisfy the remaining
inequalities. Note that the second inequality is equivalent to:
$$
\left(\frac{2}{\delta M}\right)^k\leq \frac{\lambda\delta}{4\|P(1)\|_{L^2}}.
$$
It is so sufficient to take:
$$
M\geq \sup\left(\frac{2}{\delta}, \frac{8\|P(1)\|_{L^2}}{\lambda \delta^2}\right).
$$

The third inequality is equivalent to:
$$
\left(\frac{M}{C_0^N}\right)^{k/N}\geq M^{3(1+1/N)}.
$$
For this case it is sufficient to take $M\geq \sup(1,C_0^{2N})$. Indeed,
this ensures $M^2/C_0^{2N}\geq M$ and so:
$$
\left(\frac{M}{C_0^N}\right)^{k/N}\geq M^{k/(2N)}\geq M^6,
$$
for $k\geq 12N$. But $M^6\geq M^{3(1+1/N)}$ for $M\geq 1$ and $N\geq 2$.

Therefore we can fix:
$$
M=\sup\left(1,C_0^{2N},\frac{2}{\delta}, \frac{8\|P(1)\|_{L^2}}{\lambda \delta^2} \right).
$$

The constant $\lambda$, $\delta$, and $M$ are now
fixed for the rest of the proof. The constant $\eps_0$ will be
constructed from those ones. \vskip0.2cm \noindent{\bf Step 3.
Induction:} We set:
$$
\tek=(\te-C_k)_+,
$$
with $C_k=2-\lambda(1+2^{-k})$. We consider a cut-off function in
$x$ only such that:
\begin{eqnarray*}
&& \un_{\{B_{1+2^{-k-1}}\}}\leq\eta_k\leq
\un_{\{B_{1+2^{-k}}\}},\\
&&|\nabla\eta_k|\leq C 2^k,
\end{eqnarray*}
and we denote:
$$
A_k=2\int_{-1-2^{-k}}^{0}\int_0^{\delta^k}\nint|\nabla(\eta_k\tesk)|^2\,dx\,dz\,dt+\sup_{[-1-2^{-k},1]}\nint(\eta_k\tek)^2\,dx\,dt.
$$
We want to prove simultaneously that for every $k\geq0$:
\begin{eqnarray}
\label{2} &&A_k\leq M^{-k}\\ \label{3}
 &&\eta_k\tesk \ \ \mathrm{is \
supported \ in \ } 0\leq z\leq \delta^k.
\end{eqnarray}
 \vskip0.2cm \noindent {\bf Step
4. Initial step:} We prove in this step that (\ref{2}), 
is verified for $0\leq k\leq 12N$, and that (\ref{3}) is verified
for $k=0$. We use the energy inequality (\ref{eq_energy_locale})
with cut-off function $\eta_k(x)\psi(z)$ where $\psi$ is a fixed
cut-off function in $z$ only.  Taking the mean value of
(\ref{eq_energy_locale}) in $t_1$ between $-4$ and $-2$, we find
that (\ref{2}) is verified for $0\leq k\leq 12N$ if $\eps_0$ is
taken such that:
\begin{equation}\label{alpha}
C2^{24N}(1+\Phi)\eps_0\leq M^{-12N}.
\end{equation}
We have used that $|\nabla\eta_k|^2\leq C2^{24N}$ for $0\leq k\leq 12N$. Let
us consider now  the support property (\ref{3}).
By the maximum principle, we have:
$$
\tes\leq (\te_+\un_{B_4})\ast P(z)+ b_1(x,z),
$$
in $\R^+\times B^*_4$, where $P(z)$ is the Poisson kernel. Indeed,
the right-hand side function is harmonic, positive and the trace on
the boundary is bigger that the one of $\tes$.

From step 1 we have: $b_1(x,z)\leq 2-4\lambda$. Moreover:
$$
\|\te_+\un_{B_4}\ast P(z)\|_{L^\infty(z\geq1)}\leq
C\|P(1)\|_{L^2}\sqrt{\eps_0}\leq C\sqrt{\eps_0}.
$$
Choosing $\eps_0$ small enough such that this constant is smaller
that $2\lambda$ gives:
$$
\tes\leq 2-2\lambda\qquad \mathrm{for} \ \ z\geq1, t\geq0, x\in B,
$$
so:
$$
\tes_0=(\tes-(2-2\lambda))_+\leq0 \qquad \mathrm{for} \ \ z\geq1,
t\geq0, x\in B.
$$
Hence $\eta_0\tes_0$ is supported in $0\leq z\leq \delta^0=1$.
\vskip0.2cm\noindent{\bf Step 5. Propagation of the support property
(\ref{3}):} Assume that (\ref{2}) and (\ref{3}) are verified at $k$.
We want to show that (\ref{3}) is verified at $(k+1)$. We will show
also that the following  is verified at $k$:
\begin{equation}\label{4}
\eta_k\tes_{k+1}\leq [(\eta_k\tek)\ast P(z)]\eta_k.
\end{equation}
 We consider the set $\Qks=B_{1+2^{-k}}\times[0,\delta^k]$, and we
want to control $\tesk$ on this set by harmonic functions taking
into account the contributions of the sides one by one. On
$z=\delta^k$ we have no contribution thanks to the induction
property (\ref{3}) at $k$ (the trace is equal to 0). The
contribution of the side $z=0$ can be controlled by: $\eta_k\tek\ast
P(z)$ (It has the same trace than $\tek$ on $B_{1+2^{-k-1}}$).

On each of the other side we control the contribution by:
$$
b_2((x_i-x^+)/\delta^{k},z/\delta^{k})+b_2((-x_i+x^-)/\delta^{k},z/\delta^{k}),
$$
where $x^+=(1+2^{-k})$ and $x^-=-x^+$. Indeed, $b_2$ is harmonic,
and on the side $x_i^+$ and $x_i^-$ it is bigger than 2. Finally, by
the maximum principle:
$$
\tesk\leq \sum_{i=1}^{N}\left[b_2((x_i-x^+)/\delta^{k},z/\delta^{k})
+b_2((-x_i+x^-)/\delta^{k},z/\delta^{k})\right]+(\eta_k\tek)\ast
P(z).
$$
From Step 1, for $x\in B_{1+2^{-k}}$:
\begin{eqnarray*}
&&\sum_{i=1}^{N}\left[b_2((x_i-x^+)/\delta^{k},z/\delta^{k})
+b_2((-x_i+x^-)/\delta^{k},z/\delta^{k})\right]\\
&&\qquad \leq N \Cb e^{-\frac{\pi 2^{-k}}{\delta^k}}\\
&&\qquad \leq \lambda 2^{-k-2},
\end{eqnarray*}
(thanks to Step 2). This gives (\ref{4}) since:
$$
\tes_{k+1}\leq (\tesk-\lambda 2^{-k-1})_+.
$$
More precisely, this gives:
$$
\tes_{k+1}\leq ((\eta_k\tek)\ast P(z)-\lambda 2^{-k-2})_+.
$$
 So:
$$
\eta_{k+1}\tes_{k+1}\leq ((\eta_k\tek)\ast P(z)-\lambda 2^{-k-2})_+.
$$
From the second property of Step 2, we find for $z=\delta^{k+1}$:
\begin{eqnarray*}
|(\eta_k\tek)\ast P(z)|&\leq& A_k \|P(z)\|_{L^2}\\
&\leq&\frac{M^{-k}}{\delta^{k+1}}\|P(1)\|_{L^2}\leq \lambda
2^{-k-2}.
\end{eqnarray*}
The last inequality makes use of Step 2. Therefore:
$$
\eta_{k+1}\tes_{k+1}\leq0\qquad \mathrm{on} \ \ z=\delta^{k+1}.
$$
Note, in particular, that with step 4 this gives that (\ref{3}) is
verified up to $k=12N+1$ and (\ref{4}) up to $k=12N$.
\vskip0.2cm\noindent{\bf Step 6. Propagation of Property (\ref{2}).}
We show in this step that if (\ref{3}) is true for $k-3$ and
(\ref{2}) is true for $k-3$, $k-2$ and $k-1$ then (\ref{2}) is true
for $k$.

First notice that from Step 5, (\ref{3}) is true at $k-2$, $k-1$,
and $k$.  We just need  to show that:
\begin{equation}\label{1}
A_{k}\leq C_0^k(A_{k-3})^{1+1/N}\qquad \mathrm{for} \ k\geq 12N+1,
\end{equation}
with:
\begin{equation}\label{eq_C_0}
C_0=C\frac{2^{1+2/N}}{\lambda^{2/N}}.
\end{equation}
 Indeed, the third inequality of  Step 2 gives the result.
 \vskip0.2cm\noindent{\bf Step 7. Proof of (\ref{1}):} Since
 $\eta\tes_+$ has the same trace at $z=0$ that $(\eta \te_+)^*$ and
 the latter is harmonic we have:
 $$
\int |\na(\eta\tes_+)|^2\geq
\int|\nabla(\eta\te_+)^*|^2=\int|\La^{1/2}(\eta\te_+)|^2.
 $$
Sobolev and
 Holder inequalities give:
 $$
A_{k-3}\geq
C\|\eta_{k-3}\te_{k-3}\|^2_{L^{\frac{2(N+1)}{N}}([-1-2^{-k-3},0]\times\R^N)}.
 $$
 From (\ref{4}):
 $$
\|\eta_{k-3}\tes_{k-2}\|^2_{L^{\frac{2(N+1)}{N}}}\leq
\|P(1)\|^2_{L^1}\|\eta_{k-3}\te_{k-3}\|^2_{L^{\frac{2(N+1)}{N}}}.
 $$
 So:
 \begin{eqnarray*}
A_{k-3}&\geq&C\|\eta_{k-3}\tes_{k-2}\|^2_{L^{\frac{2(N+1)}{N}}}+C\|\eta_{k-3}\te_{k-3}\|^2_{L^{\frac{2(N+1)}{N}}}\\
&\geq&
C\left(\|\eta_{k-1}\tes_{k-1}\|^2_{L^{\frac{2(N+1)}{N}}}+\|\eta_{k-1}\te_{k-1}\|^2_{L^{\frac{2(N+1)}{N}}}\right).
 \end{eqnarray*}
 Taking the mean value of (\ref{eq_energy_locale}) in $t_1$ between
 $-1-2^{-k-1}$ and  $-1-2^{-k}$, we find:
 $$
A_k\leq
C2^k(\Phi+2)\left(\int\eta_{k-1}^2\tek^2+\int\eta_{k-1}^2{\tesk}^2\right).
 $$
 Note that we have used here (\ref{3}) since $\eta_{k-1}$ is a
 cut-off function in $x$ only.
If $\te_{k}>0$ then $\te_{k-1}\geq 2^{-k}\lambda$. So:
$$
\un_{\{\te_{k}>0\}}\leq \frac{C2^{k}}{\lambda}\te_{k-1}.
$$
Therefore:
$$
A_{k}\leq \frac{C2^{k(1+2/N)}}{\lambda^{2/N}}A_{k-3}^{1+1/N}.
$$
This gives (\ref{1}).\qquad\qed

\section{The second technical lemma.}\label{section_lemme2}

We set $Q_r=B_r\times[-r,0]$ and $Q^*_r=B^*_r\times[-r,0]$

\begin{lemm}\label{lemm_2technical}
For every $\eps_1>0$, there exists a constant $\delta_1>0$ with the
following property:

For every solution $\te$ to (\ref{eq_quasi_linear}) with $v$
verifying (\ref{condition sur v}) and:
\begin{eqnarray*}
&&\tes\leq 2 \qquad\mathrm{in} \ \ Q^*_{4},\\
&& \left|\{(x,z,t)\in Q_4^*; \ \tes(x,z,t)\leq 0
\}\right|\geq\frac{|Q_4^*|}{2},
\end{eqnarray*}
we have the following implication:
$$
 \left|(x,z,t)\in
Q_4^*; \ 0<\{\tes(x,z,t)< 1 \}\right|\leq\delta_1
$$
implies:
$$
\int_{Q_{1}}(\te-1)_+^2\,dx\,dt+\int_{Q^*_{1}}({\tes}-1)_+^2\,dx\,dx\,dt\leq
C\sqrt{\eps_1}.
$$
\end{lemm}
\pf Take $\eps_1\ll1$. From the energy inequality
(\ref{eq_energy_locale}), we get:
$$
\int_{-4}^0\int_{B^*_1}|\nabla \tes_+|^2\,dx\,dz\,dt\leq C.
$$
Let:
$$
K=\frac{4\int|\nabla\tes_+|^2\,dx\,dz\,dt}{\eps_1}.
$$
Then:
\begin{equation}\label{eq_petit gradients}
\left|\{t | \int_{B^*_1}|\nabla\tes_+|^2(t)\,dx\,dz\geq
K\}\right|\leq\frac{\eps_1}{4}.
\end{equation}
For all $t\in \{t | \int_{B^*_1}|\nabla\tes_+|^2(t)\,dx\,dz\leq
K\}$, the De Giorgi lemma (see appendix) gives that:
$$
|\At||\Bt|\leq |\Ct|^{1/2}K^{1/2},
$$
where:
\begin{eqnarray*}
&&\At=\{(x,z)\in B_1^*\  | \ \tes(t,x,z)\leq0\}\\
&&\Bt=\{(x,z)\in B^*_1\  | \ \tes(t,x,z)\geq1\}\\
&&\Ct=\{(x,z)\in B^*_1\  | \  0<\tes(t,x,z)<1\}.
\end{eqnarray*}
Let us set
$$\delta_1=\eps_1^8,$$
$$I=\{t\in[-4,0]; \  |\Ct|^{1/2}\leq \eps_1^3 \ \ \mathrm{and}
 \ \ \int_{B^*_1}|\nabla \tes_+|^2(t)\,dx\,dz\leq K\}.$$
First we have, using Tchebichev inequality:
\begin{eqnarray*}
&&\left|\{t\in[-4,0]; \
|\Ct|^{1/2}\geq \eps_1^3\}\right|\\
&&\qquad \leq \frac{|\{(t,x,z) \ | \
0<\tes<1\}|}{\eps_1^6}\\
&&\qquad\leq\frac{\delta_1}{\eps_1^6}\\
&&\leq \eps_1^2\leq \eps_1/4.
\end{eqnarray*}
Hence $|[-4,0]\setminus I|\leq\eps_1/2$. Secondly we get for every
$t\in I$ such that $|\At|\geq 1/4$:
\begin{equation}\label{*}
 |\Bt|\leq
\frac{|\Ct|^{1/2} K^{1/2}}{\At}\leq 4\eps_1^{5/2}\leq \eps_1^2.
\end{equation}
 In particular:
 \begin{eqnarray*}
 \int{\tes_+}^2(t)\,dx\,dz&\leq&4(|\Bt|+|\Ct|)\\
 &\leq&8\eps_1^2.
\end{eqnarray*}
And so:
$$
\int\te_+^2(t)\,dx\leq \sqrt{K}\sqrt{\int{\tes_+}^2(t)\,dx\,dz}\leq
C\sqrt{\eps_1}.
$$
 We want to show that $|\At|>1/4$ for every $t\in I\cap [-1,0]$. First,
since $|\{(t,x,z)\ | \ \tes\leq0\}|\geq |Q_4^*|/2$, there exists
$t_0\leq-1$ such that $|{\mathcal{A}}(t_0)|\geq 1/4$. So for this
$t_0$, $\int\te_+^2(t_0)\,dx\leq C\sqrt{\eps_1}$. Using the energy
inequality (\ref{eq_energy_locale}), we have for every $t\geq t_0$:
$$
\int\te_+^2(t)\,dx\leq \int\te_+^2(t_0)\,dx+C(t-t_0).
$$
So for $t-t_0\leq\delta^*=1/(64C)$ we have:
\begin{eqnarray*}
&& \int\te_+^2(t)\,dx\leq \frac{1}{64}.
\end{eqnarray*}
(Note that $\delta^*$ do not depend on $\eps_1$. Hence we can
suppose $\eps_1\ll \delta^*$.) We have:
\begin{eqnarray*}
\tes_+(z)&\leq&\te_++\int_0^z\partial_z\tes_+\,dz\\
&\leq&\te_++\sqrt{z}\left(\int|\partial_z\tes_+|^2\,dz\right)^{1/2}.
\end{eqnarray*}
So, for $t-t_0\leq\delta^*$, $t\in I$ and $z\leq \eps_1^2$ we have:
$$
\tes_+(t,x,z)\leq
\te_+(t,x)+\left(\eps_1^2\int|\partial_z\tes_+|^2\,dz\right)^{1/2}.
$$
The integral in $x$ of the right hand side term is less than
$1/8+\sqrt{\eps_1}\leq 1/4$. So by Tchebichev:
$$
|\{z\leq\eps_1^2, \ x\in B_1, \ \tes_+(t)\geq1\}|\leq
\frac{\eps^2_1}{4}.
$$
Since $|\Ct|\leq\eps_1^6$, this gives
$$
|\At|\geq \eps_1^2(1-1/4)-\eps_1^6\geq \eps_1^2/2.
$$
Then (\ref{*}) gives:
$$
|\Bt|\leq 2\sqrt{\eps_1},
$$
and:
$$
|\At|\geq1-2\sqrt{\eps_1}-\eps_1^6\geq1/4.
$$
Hence, for every $t\in[t_0,t_0+\delta^*]\cap I$ we have:
$|\At|\geq1/4$. On $[t_0+\delta^*/2,t_0\delta^*]$ there exists
$t_1\in I$ ($\delta^*\geq\eps_1/4$). And so, we can construct an
increasing sequence $t_n$, $0\geq t_n\geq t_0+n\delta^*/2$ such that
$|\At|\geq1/4$ on $[t_n,t_n+\delta^*]\cap I\supset[t_n,t_{n+1}]\cap
I$. Finally on $I\cap[-1,0]$ we have $|\At|\geq1/4$. This gives from
(\ref{*}) that for every $t\in I\cap[-1,0]$: $|\Bt|\leq \eps_1/16$.
Hence:
$$
|\{\tes\geq1\}|\leq \eps_1/16+\eps_1/2\leq\eps_1.
$$
Since $(\tes-1)_+\leq 1$, this gives that:
$$
\int_{Q^*_1}(\tes-1)^2_+\,dx\,dz\,dt\leq \eps_1.
$$
We have for every $t,x$ fixed:
$$
\te-\tes(z)=\int_0^z\partial_z\tes\,dz.
$$
So:
\begin{eqnarray*}
(\te-1)_+^2&\leq&
2\left(({\tes}(z)-1)_+^2+\left(\int_0^z|\nabla\tes|\,dz\right)^2\right)\\
&\leq&\frac{2}{\sqrt{\eps_1}}\int_0^{\sqrt{\eps_1}}({\tes}-1)_+^2\,dz+2\sqrt{\eps_1}\int_0^{\sqrt{\eps_1}}|\nabla
\tes|^2\,dz.
\end{eqnarray*}
Therefore:
$$
\int_{Q_1}(\te-1)^2_+\,dx\,ds\leq C\sqrt{\eps_1}.
$$
 \qed

\section{Oscillation lemma}\label{section_oscillation lemma}
 This section is dedicated to the proof of the following
 proposition:
 \begin{prop}\label{lemm_oscillations}
There exits $\lambda^*>0$ such that for every solution $\te$ of
(\ref{eq_quasi_linear}) with $v$ verifying (\ref{condition sur v}),
if:
\begin{eqnarray*}
&&\tes\leq2\qquad \mathrm{in} \ \ Q^*_1\\
&&|\{(t,x,z)\in Q^*_1; \ \tes\leq0\}|\geq\frac{1}{2},
\end{eqnarray*}
then:
$$
\tes\leq 2-\lambda^*\qquad\mathrm{in} \ \  Q^*_{1/16}.
$$
\end{prop}
\pf For every $k\in\N$, $k\leq K_+=E(1/\delta_1+1)$ (where
$\delta_1$ is defined in Lemma \ref{lemm_2technical} for $\eps_1$
such that $4C\sqrt{\eps_1}\leq\eps_0$, $\eps_0$ defined in Lemma
\ref{lemm_L^2_L^infty}), we define:
$$
\tbk=2(\tbku-1)\qquad \mathrm{with}\qquad \overline{\te}_0=\te.
$$
So we have: $\tbk=2^k(\te-2)+2$. Note that for every $k$, $\tbk$
verifies (\ref{eq_quasi_linear}), $\tbk\leq2$ and $|\{(t,x,z)\in
Q^*_1 \ | \ \tbk\leq0\}|\geq\frac{1}{2}$. Assume that for all those
$k$, $|\{0<\tbsk<1\}|\geq\delta_1$. Then, for every $k$:
$$
|\{\tbsk<0\}|=|\{\tbsku<1\}|\geq|\{\tbsku<0\}|+\delta_1.
$$
Hence:
$$
|\{\overline{\te}^*_{K_+}\leq 0\}|\geq1,
$$
and $\overline{\te}^*_{K_+}<0$ almost everywhere, which means:
$2^{K_+}(\tes-2)+2<0$ or
$$
\tes<2-2^{-K_+}.
$$
And in this case we are done.

Else, there exists $0\leq k_0\leq K_+$ such that:
$|\{0<\overline{\te}^*_{k_0}<1\}|\leq\delta_1$. From Lemma
\ref{lemm_2technical} and Lemma \ref{lemm_L^2_L^infty} (applied on
$\overline{\te}_{k_0+1}$) we get $(\overline{\te}_{k_0+1})_+\leq
2-\lambda$ which means:
$$
\te\leq 2-2^{-(k_0+1)}\lambda\leq 2-2^{-K_+}\lambda,
$$
in $Q_{1/8}$.

Consider the function $b_3$ defined by:
\begin{eqnarray*}
&&\Delta b_3=0\qquad \mathrm{in} \ \ B^*_{1/8}, \\
&&b_3=2\qquad \mathrm{on \ the \ sides \ of \ the \ cube \ except \
 for \ }z=0\\
 && b_3=2-2^{-K_+}\inf(\lambda,1) \qquad \mathrm{on} \ z=0.
\end{eqnarray*}
We have $b_3<2-\lambda^*$ in $B^*_{1/16}$. And from the maximum
principle we get $\te^*\leq b_3$.\qquad\qed

\section{Proof of Theorem \ref{theo L^infty to C^alpha}.}\label{section_theorem}

We fix $t_0>0$ and consider $t\in [t_0,\infty[\times\R^N$. We
define:
$$
F_0(s,y)=\te(t+st_0/4,x+t_0/4(y-x_0(s))),
$$
where $x_0(s)$ is solution to:
\begin{eqnarray*}
&&\dot{x}_0(s)=\frac{1}{|B_4|}\int_{x_0(s)+B_4}v(t+st_0/4,x+yt_0/4)\,dy\\
&&x_0(0)=0.
\end{eqnarray*}
Note that $x_0(s)$ is uniquely defined from Cauchy Lipschitz
theorem. We set:
\begin{eqnarray*}
&&\tilde{\te}^*_0(s,y)=\frac{4}{\sup_{Q^*_4} F^*_0-\inf_{Q^*_4}
F^*_0}\left(F^*_0-\frac{\sup_{Q^*_4}
F^*_0+\inf_{Q^*_4} F^*_0}{2}\right).\\
&& v_0(s,y)=v(t+st_0/4,x+t_0/4(y-x_0(s)))-\dot{x}_0(s),
\end{eqnarray*}
and then for every $k>0$:
\begin{eqnarray*}
&&F_k(s,y)=F_{k-1}(\tilde{\mu}s,\tilde{\mu}(y-x_k(s))),\\
&&\tilde{\te}^*_k(s,y)=\frac{4}{\sup_{Q^*_4} F^*_k-\inf_{Q^*_4}
F^*_k}\left(F^*_k-\frac{\sup_{Q^*_4}
F^*_k+\inf_{Q^*_4} F^*_k}{2}\right),\\
&&\dot{x}_k(s)=\frac{1}{|B_4|}\int_{x_k(s)+B_4}v_{k-1}(\tilde{\mu}s,\tilde{\mu}y)\,dy\\
&&x_k(0)=0\\
&&v_k(s,y)=v_{k-1}(\tilde{\mu}s,\tilde{\mu}(y-x_k(s)))-\dot{x}_k(s),
\end{eqnarray*}
where $\tilde{\mu}$  will be chosen later. We divide the proof in
several steps. \vskip0.2cm\noindent{\bf Step 1.} For k=0,
$\tilde{\te}_0$ is solution to (\ref{eq:quasi}) in $[-4,0]\times
\R^N$, $\|v_0\|_{BMO}=\|v\|_{BMO}$, $\int v_0(s)\,dy=0$ for every
$s$ and $|\tilde{\te}_0|\leq 2$. Assume that it is true at $k-1$.
Then:
$$
\partial_s
F_k=\tilde{\mu}\partial_s\tilde{\te}_{k-1}()-\tilde{\mu}\dot{x}_k(s)\cdot\nabla\tilde{\te}_{k-1.}
$$
So $\tilde{\te}_k$ is solution of (\ref{eq:quasi}) and
$|\tilde{\te}_k|\leq 2$. By construction, for every $s$ we have
$\int_{B_4}v_k(s,y)\,dy=0$ and
$\|v_k\|_{BMO}=\|v_{k-1}\|_{BMO}=\|v\|_{BMO}$. Moreover we have:
\begin{eqnarray*}
|\dot{x}_k(s)|&\leq& \int_{B_4}v_{k-1}(\tilde{\mu}(y-x_k(s)))\,dy\\
&\leq&C\|v_{k-1}(\tilde{\mu}y)\|_{L^p}\\
&\leq&C\tilde{\mu}^{-N/p}\|v_{k-1}\|_{L^p}\\
&\leq&C_p\tilde{\mu}^{-N/p}\|v_{k-1}\|_{BMO}.
\end{eqnarray*}
So, for $0\leq s\leq1$, $y\in B_4$ and $p>N$:
$$
|\tilde{\mu}(y-x_k(s))|\leq
4\tilde{\mu}(1+C_p\tilde{\mu}^{-N/p})\leq C\tilde{\mu}^{1-N/p}.
$$
For $\tilde{\mu}$ small enough this is smaller than 1.
\vskip0.2cm\noindent{\bf Step 2.} For every $k$
 we can use the oscillation lemma. If $|\{\tilde{\te}^*_k\leq0\}|\geq
 \frac{1}{2}|Q_4^*|$ then we have $\tilde{\te}^*_k\leq
 2-\lambda^*$. Else we have $|\{-\tilde{\te}^*_k\leq0\}|\geq
 \frac{1}{2}|Q_4^*|$ and applying the oscillation lemma on
 $-\tilde{\te}^*_k$ gives $\tilde{\te}^*_k\geq-2+\lambda^*$.
 In both cases this gives:
 $$
|\sup\tilde{\te}^*_k-\inf\tilde{\te}^*_k|\leq2-\lambda^*.
 $$
 and so:
$$
|\sup_{Q_1^*}F^*_k-\inf_{Q^*_1}F^*_k|\leq(1-\lambda^*/2)^k|\sup_{Q_1^*}F^*_0-\inf_{Q^*_1}F^*_0|.
$$
\vskip0.2cm\noindent{\bf Step 3.} For $s\leq \tilde{\mu}^{2n}$:
$$
\sum_{k=0}^{n}\tilde{\mu}^{n-k}x_k(s)\leq
\tilde{\mu}^{2n}\sum_{k=0}^n
\frac{\tilde{\mu}^{n-k}}{\tilde{\mu}^{-N/p}}\leq\frac{\tilde{\mu}^n}{2},
$$
for $\tilde{\mu}$ small enough. So
$$
\left|\sup_{[-\tilde{\mu}^{2n},0]\times
B^*_{\tilde{\mu}^{n}/2}}\tes-\inf_{[-\tilde{\mu}^{2n},0]\times
B^*_{\tilde{\mu}^{n}/2}}\tes\right|\leq (1-\lambda^*/2)^{n}.
$$
This gives that  $\tes$ is $C^\alpha$ at $(t,x,0)$, and so $\te$ is
$C^\alpha$ at $(t,x)$.\qquad\qed
\appendix
\section{Proof of the De Giorgi isoperimetric lemma.}

Let $\omega\in H^1([-1,1]^{N+1})$. We denote:
\begin{eqnarray*}
&&{\mathcal{A}}=\{x; \ \omega(x)\leq0)\}\\
&&{\mathcal{B}}=\{x; \  \omega(x)\geq1)\}\\
&&{\mathcal{C}}=\{x; \  0<\omega(x)<1)\},
\end{eqnarray*}
and
$$
\chi=\un_{\{y_1+s(y_1-y_2)/|y_1-y_2|\in{\mathcal{C}}\}}.
$$
 We have:
\begin{eqnarray*} &&|\mathcal{A}||\mathcal{B}|\leq
\int_{\mathcal{A}}\int_{\mathcal{B}}(\omega(y_1)-\omega(y_2))\,dy_1\,dy_2\\
&&\qquad
=\int_{\mathcal{A}}\int_{\mathcal{B}}\int_0^{|y_1-y_2|}\nabla\omega(y_1+s\frac{y_1-y_2}{|y_1-y_2|})\cdot\frac{y_1-y_2}{|y_1-y_2|}\,ds\,dy_1\,dy_2\\
&&\qquad=\int_{\mathcal{A}}\int_{\mathcal{B}}\int_0^{|y_1-y_2|}\chi\nabla\omega(y_1+s\frac{y_1-y_2}{|y_1-y_2|})\cdot\frac{y_1-y_2}{|y_1-y_2|}\,ds\,dy_1\,dy_2\\
&&\qquad
\leq\int_{\mathcal{A}}\int_{\mathcal{B}}\int_0^\infty\chi\left|\nabla\omega(y_1+s\frac{y_1-y_2}{|y_1-y_2|})\right|\,ds\,dy_1\,dy_2\\
&&\qquad
\leq\int_{B_1}\int_{B_1}\int_0^\infty\chi\left|\nabla\omega(y_1+s\frac{y_1-y_2}{|y_1-y_2|})\right|\,ds\,dy_1\,dy_2\\
&&\qquad
\leq\int_{S_{N-1}}\int_{B_1}\int_0^\infty\frac{\left|\nabla\omega(y_1+s\nu)\right|}{s^{N-1}}\un_{\{(y_1+s\nu)\in\mathcal{C}\}}s^{N-1}\,d\nu\,ds\,dy_1\\
&&\qquad
\leq C\int_{B_1}\int_{B_1}\frac{\left|\nabla\omega(y_1+y_2)\right|}{|y_2|^{N-1}}\un_{\{y_1+y_2\in\mathcal{C}\}}\,dy_2\,dy_1\\
&&\qquad \leq C\|\nabla\omega\|_{L^2}|\mathcal{C}|^{1/2}.
\end{eqnarray*}
\qed

\section{Higher regularity}

We give the proof of the following theorem.
\begin{theo}
Let $\te$ be a solution of the quasi-geostrophic equation
(\ref{eq:quasi}), (\ref{eq:uquasi}) satisfying the regularity
properties of Theorem \ref{theo quasi}:
\begin{eqnarray*}
\te&\in&L^\infty(0,\infty;L^2)\cap L^2(0,\infty;H^{1/2})\\
&&\cap L^\infty([t_0,\infty[\times\R^N)\cap
C^\alpha([t_0,\infty[\times\R^N),
\end{eqnarray*}
for every $t_0>0$. Then $\te$ belongs to
$C^{1,\beta}([t_0,\infty[\times\R^N)$ for every $\beta<1$ and
$t_0>0$ and is therefore a classical solution.
\end{theo}
\noindent{\bf Proof:} We want to show the regularity at a fixed
point $y_0=(t_0,x_0)\in ]0,\infty[\times\R^N\subset \R^m$ where
$m=N+1$. Note that Changing $\te(t,x)$ by
$\te(t,x-u(t_0,x_0)t)-\te(t_0,x_0)$ if necessary, we can assume
without loss of generality that $\te(y_0)=0$ and $u(y_0)=0$. The
fundamental solution of:
$$
\dt \te+\La\te=0
$$
is the Poisson kernel:
$$
P(t,x)=\frac{Ct}{(|x^2|+t^2)^{\frac{N+1}{2}}},
$$
a homogeneous function of order $-N$ if extended for $t$ negative.
the solution $\te$ of (\ref{eq:quasi}) 
can be represented  as the sum of two terms.
\begin{equation}\label{eq_2terms}
\te(t,x)=P(t,\cdot)\ast\te_0-g(t,x),
\end{equation}
where:
\begin{eqnarray*}
g(t,x)&=&\int_0^t\nint
P(t-t_1,x-x_1)\Div(u(t_1,x_1)\te(t_1,x_1))\,dt_1\,dx_1\\
&=&\int_0^\infty\int_{\R^N}\nabla_x\tP(y-y_1)\cdot
u(y_1)\te(y_1)\,dy_1.
\end{eqnarray*}
In the last inequality, we denoted $y=(t,x)$, $\tP$ the extension of
$P$ for negative $t$ with  value 0, and  we passed the divergence on
$\tP$, which becomes a singular integral.
 The first term in (\ref{eq_2terms}) is smooth for $t>0$ and depends only
 on the initial data. We
focus on the second one $g(y)$. We fix $e\in {\mathbf{S}}_m$, and
estimate
 $g(y_0+he)-g(y_0)$ for $h>0$ in the standard way. We split the
 integral:
\begin{equation}\label{eq_delta g}
g(y_0)-g(y_0+he)
=\int_0^\infty\int_{\R^N}Q_0(y_0-y_1,he)u(y_1)\te(y_1)\,dy_1
 \end{equation}
 where:
 $$
Q_0(y,he)=\nabla_x\tP(y)-\nabla_x\tP(y+he),
 $$
into two parts, one on the ball $B_{10h}$ centered to $y_0$ and
radius $10h$, 
and the second on the complement.
The first part has no cancelation  so we separate the integrals:
\begin{eqnarray*}
&&\qquad\int_{B_{10h}}\un_{\{t_1\geq0\}}
[\nabla_x\tP(y_0-y_1)-\nabla_x\tP(y_0+he-y_1)]u(y_1)\te(y_1)\,dy_1\\
&&=\int_{B_{10h}}\un_{\{t_1\geq0\}}\nabla_x\tP(y_0-y_1)u(y_1)\te(y_1)\,dy_1\\
&&\qquad-\int_{B_{10h}}\un_{\{t\geq0
\}}\nabla_x\tP(y_0+he-y_1)u(y_1)\te(y_1)\,dy_1.
\end{eqnarray*}
If $\te$ is $C^\alpha$, $\alpha>0$, from the Riesz transform $u$ is
also $C^\alpha$, and since $\te(y_0)=u(y_0)=0$, we have:
\begin{equation}\label{eq_utca}
|u(y_1)\te(y_1)|\leq \inf(|y_1-y_0|^{2\alpha},C).
\end{equation}
So the first integral is convergent and bounded by $Ch^{2\alpha}$.
To deal with the second one, notice that $\nabla_x \tP$ have mean
value zero on any slice $t=C$ of $B_{10h}$, so we can add and
substract $\te(y_0+he)u(y_0+he)$. We have:
\begin{eqnarray*}
&&|\te(y_1)u(y_1)-\te(y_0+he)u(y_0+he)|\leq C
h^\alpha|y_0+he-y_1|^\alpha,
\end{eqnarray*}
where we have used again that $u(y_0)=\te(y_0)=0$. Hence the
integral is also convergent and bounded by $Ch^{2\alpha}$. This
gives that the contribution of $B_{10h}$ on (\ref{eq_delta g}) is
smaller that $C h^{2\alpha}$.
%

Outside of a neighborhood of size $10 h$  we use the cancelation of
$\nabla_x \tP$. Up to Lipschitz regularity we just do:
\begin{eqnarray*}
&& |\nabla_x[\tP(y_1-y_0)-\tP(y_1+he-y_0)]|\\
&&\qquad\leq \frac{h}{|y_1-y_0|^{m+1}},
\end{eqnarray*}
and integrate against $|u\te|$ which verifies (\ref{eq_utca}). This
gives the bound:
\begin{eqnarray*}
&&\int_{|y_1-y_0|\geq10h}\frac{h}{|y_1-y_0|^{m+1-2\alpha}}\,dy_1
 \leq Ch^{2\alpha},
\end{eqnarray*}
provided that $2\alpha<1$. Altogether, this gives that if $\te\in
C^\alpha$ with $2\alpha<1$, then
$$
|g(y_0)-g(y_0+he)|\leq Ch^{2\alpha}.
$$
Bootstrapping the argument gives that $\te$ is $C^\alpha$ for any
$\alpha<1$.

 To go beyond Lipschitz we consider a second order
increment quotient:
$$
Q_1(y,he)=|\nabla[\tP(y+he)+\tP(y-he)-2\tP(y)]|.
$$
We have:
$$
g(y_0+he)+g(y_0-he)-2g(y_0)=\int_{\R^m}\un_{\{t_1\geq0\}}Q_1(y_0-y_1,he)u(y_1)\te(y_1)\,dy_1.
$$
Note that $Q_1(y,he)=Q_0(y,he)-Q_0(y-he,he)$, so for $|y|<20h$, the
local estimate of the previous argument together with the $C^\alpha$
property of $\te$ and $u$ gives:
$$
\int_{B_{20h}}|Q_1(y_0-y_1)u(y_1)\te(y_1)|\,dy_1\leq Ch^{2\alpha}.
$$
For $|y|> 20 h$ and $y$ not in the strip
${\mathcal{T}}_h=[t_0-h,t_0+h]\times\R^N$, we have:
$$
|Q_1(y_0-y_1,he)|\leq C\frac{h^2}{|y_0-y_1|^{m+2}}.
$$
and the corresponding integral:
\begin{eqnarray*}
&&\int_{|y_0-y_1|\geq 20
h}\un_{\{y_1\notin{\mathcal{T}}_h\}}|Q_1(y_0-y_1)u(y_1)\te(y_1)|\,dy_1\\
&&\qquad\leq C\int_{|y|\geq20h}
\frac{h^2}{|y|^{m+2}}(|y|^{2\alpha}\wedge1)\,dy\\
&&\qquad\leq C h^{2\alpha},
\end{eqnarray*}
whenever $2\alpha<2$. It remains to control the contribution of the
strip ${\mathcal{T}}_h\setminus B_{20h}$. The estimate on $Q_0$
gives that on this strip:
$$
|Q_1(y_1-y_0,he)|\leq C\frac{h}{|y_1-y_0|^{N+2}}\leq
C\frac{h}{|x_1-x_0|^{N+2}}.
$$
Not that on ${\mathcal{T}}_h\setminus B_{20h}$ we have
$|x_1-x_0|\geq h$. So the contribution of this strip is bounded by:
\begin{eqnarray*}
&&\int_{t_0-h}^{t_0+h}\int_{|x_1-x_0|\geq
h}\frac{h}{|x_1-x_0|^{N+2-2\alpha}}\,dx_1\,dt_1\\
&&\leq C\frac{h^{2\alpha}}{h}\int_{t_0-h}^{t_0+h}\,dt_1\\
&&\leq C h^{2\alpha},
\end{eqnarray*}
whenever $2\alpha<2$.
 That goes all the way to $C^{1,\beta}$ for
every $\beta<1$.\qquad\qed

\section{existence of solutions to (\ref{eq_quasi_linear})}

In this appendix we sketch the existence theory of approximate
solution of the equation (\ref{eq_quasi_linear}) satisfying the
truncated energy inequalities in the hypothesis of Theorem \ref{theo
L^2 L^infty}. We start by restricting the problem to
$B_1\times[0,\infty]$ and adding an artificial diffusion term $\eps
\Delta$. We will use the eigenfunctions $\sigma_k$ and eigenvalues
$\lambda_k^2$ of the Laplacian in $B_1$, that is:
$$
\Delta \sigma_k+\lambda^2_k\sigma_k=0.
$$
Note that $\sigma_k^*(x,z)=\sigma_k(x)e^{-\lambda_kz}$ is the
harmonic extension of $\sigma_k$ for the semi-infinite cylinder
$Q_1=B_1\times[0,\infty]$ with data 0 in the lateral boundary, and:
$$
\lambda_k\sigma_k(x)=\partial_\nu\sigma_k^*(x,0),
$$
where $\partial_\nu$ is the normal derivative. Also:
$$
\int_{Q_1} \lambda_k\sigma_k^2\,dx\,dz=\int_{\partial
Q_1}\sigma_k^*\partial_\nu\sigma_k^*\,dx=\int_{Q_1}|\nabla\sigma_k^*|^2\,dx\,dz,
$$
and thus formula is also correct for any series
$$
g(x)=\sum f_k\sigma_k(x),
$$
provided that $\sum f_k^2\lambda_k$ converges, i.e., $g\in
H^{1/2}(B_1)$.

We want to solve then  in $[0,\infty[\times B_1$ the equation:
\begin{equation}\label{eq_eps}
\dt\te+\Div (v\te)=\eps\Delta\te-(-\Delta^{1/2})\te,
\end{equation}
where $-\Delta^{1/2}\te$ is understood as the operator that maps
$\sigma_k$ to $\lambda_k\sigma_k=\partial_\nu\sigma_k^*$.

For, say, $v$ bounded and divergence free, this is straightforward
using Galerkin method: Let us restrict (\ref{eq_eps}) to $\sigma_k$,
with $1\leq k\leq k_0$, i.e. we seek a function:
$$
\te=\te_{\eps,k_0}=\sum_1^{k_0}f_k(t)\sigma_k(x)
$$
that is a solution of the equation when tested against $\sigma_k$,
$1\leq k\leq k_0$. The functions $f_k$ are solutions to the
following system of ODEs:
$$
f'_k(t)=-[\eps\lambda_k^2+\lambda_k]f_k(t)+\sum_{l=1}^{k_0}
a_{kl}f_l(t),\qquad 1\leq k\leq k_0,
$$
with initial value:
$$
f_k(0)=\int_{B_1}\te_0(x)\sigma_k(x)\,dx,
$$
where:
$$
a_{kl}=\int_{B_1}v(t,x)\cdot\nabla\sigma_k(x)\sigma_l(x)\,dx.
$$
Note that, since $v$ is divergence free, the matrix $a_{kl}$ is
antisymmetric. This leads to the estimate:
\begin{eqnarray*}
&&\sum_{k=1}^{k_0}
f_k^2(t_2)+\int_{t_1}^{t_2}\sum_{k=1}^{k_0}(\eps\lambda_k^2+\lambda_k)f_k^2(s)\,ds\\
&&\qquad\qquad=\sum_{k=1}^{k_0} f_k^2(t_1).
\end{eqnarray*}
In particular $\te_{\eps,k_0}$ satisfies the energy inequality:
\begin{eqnarray*}
&& \|\te_{\eps,k_0}(t_2)\|^2_{L^2(B_1)}+\int_{t_1}^{t_2}
\left(\|\te_{\eps,k_0}(s)\|^2_{\dot{H}^{1/2}(B_1)}+\eps\|\te_{\eps,k_0}(s)\|^2_{\dot{H}^{1}(B_1)}\right)\,ds\\
&&\qquad\qquad \leq \|\te_{\eps,k_0}(t_1)\|^2_{L^2(B_1)}.
\end{eqnarray*}
Notice also that what we call $H^{1/2}(B_1)$ corresponds to the
extension of $\te$ to the half cylinder, and such:
$$
\|\te\|_{\dot{H}^{1/2}(B_1)}\geq \|\te\|_{\dot{H}^{1/2}(\R^N)}.
$$
We now pass to the limit in $k_0$ and denote $\te_\eps$ the limit.
If we test $\te_{\eps,k_0}$ with a function $\gamma\in
L^\infty(0,T;L^2(B_1))\cap L^2(0,T;H^1(B_1))$, there is no problem
in passing to the limit in the term:
$$
\int_{t_1}^{t_2}\int_{B_1}(\nabla\gamma)v\te_{\eps,k_0}\,dx\,ds,
$$
since $\te_{\eps,k_0}$ converges strongly in $L^2([0,T]\times B_1)$.
In particular, for
$\gamma=(\te_{\eps}-\lambda)_+=\te_{\eps,\lambda}$ the term
converges to:
$$
\int_{t_1}^{t_2}\int_{B_1}\nabla[\te_{\eps,\lambda}]^2 v\,dx\,ds=0,
$$
provided that $v$ is divergence free. This leads to the following
corollaries:
\begin{coro}
The function $\te_\eps$ satisfies the hypothesis of Lemma
\ref{lemm_L^2_L^infty} independently of $\eps$, and therefore:
$$
\|\te_\eps(T)\|_{L^\infty(B_1)}\leq\frac{C}{T^{N/2}}\|\te_\eps(0)\|_{L^2(B_1)}.
$$
\end{coro}

\begin{coro}
The same theorem is true for $v\in L^2([0,T]\times B_1)$
independently of the $L^2$ norm of $v$.
\end{coro}
\pf We approximate $v$ by a mollification $v_\delta$. \qquad \qed

\begin{coro}
For $\te_0$ prescribed in $L^2(\R^N)$, the same result is true in
$[0,T]\times\R^N$.
\end{coro}
\pf We may rescale the previous theorems to the ball of radius $M$
by applying them to $\overline\te(t,x)=M^{N/2}\te(Mt,Mx)$. This
change preserves the $L^2$ norm, and so we get:
$$
\sup_{B_M} M^{N/2}\te(s,y)\leq \frac{C}{(s/M)^{N/2}},
$$
or
$$
\sup_{B_M}\te(s,y)\leq \frac{C}{(s)^{N/2}},
$$
provided that $v\in L^2(B_M)$ is divergence free. Then letting  $M$
go to infinity gives the result.\qquad \qed \vskip0.3cm

\noindent{\bf Final remark:} Since all the estimates are independent
of $\eps$ we may let $\eps$ go to zero for the limit to be  weak
solution of the limiting equation, and satisfying the truncated
energy inequalities.

Note also than the same approach can be taken for higher regularity.
Indeed, the proof of higher regularity depends only on the truncated
and localized energy inequality that is also satisfied by the
$\eps$-problem. We may then pass to the limit in $\eps$ and find a
classical solution of the limiting problem.

\vskip0.2cm \noindent{\bf Acknowledgments:} Both authors were
supported in part by NSF Grants.

\bibliography{biblio}

\end{document}